\documentclass[12pt]{article}
\input epsf
\usepackage{amsmath}

\usepackage{amssymb}
\usepackage{theorem}

\usepackage{enumerate}

\usepackage{mathrsfs}

\usepackage{geometry}

\usepackage{dsfont}

\numberwithin{equation}{section}

\newtheorem{thm}{Theorem}[section]
\newtheorem{lemma}[thm]{Lemma}
\newtheorem{prop}[thm]{Proposition}

{\theorembodyfont{\rmfamily}

\newtheorem{rmk}[thm]{Remark}
}

\newcommand{\qed}{\hfill \mbox{\raggedright \rule{.07in}{.1in}}}

\newcommand{\R}{{\mathbb R}}

\newcommand{\N}{{\mathbb N}}

\def\R{\ensuremath{\mathbb R}}
\def\N{\ensuremath{\mathbb N}}

\def\l{{\rm Leb}}

\newcommand{\dif}{\mathrm{d}}

\DeclareMathOperator*{\esssup}{ess\;sup}
\DeclareMathOperator*{\essinf}{ess\;inf}

\def\dist{\ensuremath{\text{dist}}}

\def\eps {\varepsilon}

\begin{document}

\date{\today}

\title{Central limit theorems for the shrinking target problem.}

\author{Nicolai Haydn \footnote{Mathematics Department, University of Southern California, Los Angeles, 90089-1113, email: nhaydn@math.usc.edu}\\ Matthew Nicol \footnote{Department of Mathematics, University of Houston, Houston, TX 77204-3008, 
email: nicol@math.uh.edu},  Sandro Vaienti\footnote{UMR-7332 Centre de Physique Th\'{e}orique, CNRS, Universit\'{e}
d'Aix-Marseille, Universit\'{e} du Sud, Toulon-Var and FRUMAM,
F\'{e}d\'{e}ration de Recherche des Unit\'{e}s des Math\'{e}matiques de Marseille,
CPT Luminy, Case 907, F-13288 Marseille CEDEX 9, 
email: vaienti@cpt.univ-mrs.fr}, Licheng Zhang \footnote{Department of Mathematics, University of Houston, Houston, TX 77204-3008, 
email: xiyao.fei@gmail.com}\\}
\maketitle

\begin{abstract}
   Suppose  $B_i:=
  B(p,r_i)$ are nested balls of radius $r_i$ about a point $p$ in a dynamical system $(T,X,\mu)$.
  The question of
  whether $T^i x\in B_i$ infinitely often (~i.~o.) for $\mu$ a.e.\ $x$ is often
  called the shrinking target problem.  In many dynamical settings it has been shown that if $E_n:=\sum_{i=1}^n \mu (B_i)$
  diverges then there is a quantitative rate of entry and $\lim_{n\to \infty} \frac{1}{E_n} \sum_{j=1}^{n} 1_{B_i} (T^i x) \to 1$
  for $\mu$~a.e. $x\in X$. This is a self-norming type of strong law of large numbers. We establish self-norming central
  limit theorems (CLT) of the form $\lim_{ n\to \infty} \frac{1}{a_n} \sum_{i=1}^{n} [1_{B_i} (T^i x)-\mu(B_i)] \to N(0,1)$ (in distribution)
  for a variety of hyperbolic and non-uniformly hyperbolic dynamical systems, the normalization constants are
  $a^2_n \sim E [\sum_{i=1}^n 1_{B_i} (T^i x)-\mu(B_i)]^2$. Dynamical systems to which our results apply include
   smooth expanding maps of the interval, Rychlik type maps, Gibbs-Markov maps, rational maps  and,  in higher dimensions,
   piecewise expanding maps. For such central limit theorems the main difficulty is to prove that the non-stationary variance has a limit
   in probability.  \end{abstract}

\maketitle

\section{Introduction}

Suppose $(T,X,\mu)$ is an ergodic dynamical system and $B_i (p)$ is a nested sequence of balls about a point $p \in X$. Recently 
there have been many papers concerning the behavior of the almost sure limit of the normalized  sum $\frac{1}{E_n} \sum_{i=1}^{n} 1_{B_i (p)}(x)$ where
 $E_n := \sum_{i=1}^n \mu (B_i (p))$ diverges~\cite{Chernov_Kleinbock,Kim, Galatolo,Gouezel, GNO, HNPV,Kessebohmer}. 
 If the limit is known to exist almost surely then  
 $\{B_i (p)\}$ is said to satisfy the Strong Borel Cantelli property. Many of the references we mentioned consider more general sequences of sets than nested 
 balls. The study of hitting time  statistics to a sequence of nested balls is sometimes called the shrinking target problem. In this paper we 
 study self-norming central limit theorems for the shrinking target problem, namely the distribution limit of 
$\frac{1}{a_n}  \sum_{i=1}^{n} [1_{B_i (p)}- E_n]$ where $a_n$ is a sequence of norming constants. For reasons of exposition we focus on
the case where $\mu (B_i (p) )=\frac{1}{i}$, a critical case, where $E_n=\log n$.
Our results extend (with obvious modifications to the norming sequences) to balls satisfying
$\frac{C_1}{i^{\gamma_1} }\le \mu (B_i (p)) \le \frac{C_2}{i^{\gamma_2}}$ where $C_1,C_2$ are positive constants and $0<\gamma_2\le \gamma_1 \le 1$. 
The main difficulty is to establish that the non-stationary variance has a limit in probability.  Our results are limited to non-uniformly expanding systems i.e. those without a contracting direction and are based upon the  Gordin~\cite{Gordin69}
martingale approximation approach (see also~\cite{Liverani96}).

More generally, this paper is also an attempt to study the statistics of non-stationary stochastic processes arising as observations (which perhaps change over time)
on an underlying dynamical system (which may change over time).  Conze and Raugi~\cite{Conze_Raugi}
studied similar  problems for sequential expanding dynamical systems. Somewhat related  results were obtained by N\'andori, Sz\'asz and Varj\'u~\cite{NSV}
who obtained central limit theorems in  the setting in which a fixed observation $\phi: X\to \R$ was considered on a space on which a sequence of different 
transformations acted $T_i : X\to X$ act, preserving a common invariant measure $\mu$. The main difficulty in~\cite{NSV} was also controlling the variance,
 but the setting in which the underlying maps change but the observation is fixed is simpler in some respects and more difficult in others. 
 
 We obtain fairly complete results in the case in which the transfer operator with respect to the invariant measure is quasicompact  in the bounded variation norm. These results 
 are contained in  Proposition~\ref{BV} and Theorem~\ref{piecewise}. For systems in which  the  transfer operator is quasicompact in a H\"older or
 Lipschitz space we show that under the assumption we call (SP) (derived from a Gal-Koksma lemma as formulated by  Sprindzuk~\cite{Sprindzuk}) or a form of short returns assumption  called Assumption C
 we have a  central limit theorem (~Theorem~\ref{Lip}).  Assumption C and the SP property have been shown to hold for generic points in a variety of non-uniformly expanding systems~\cite{Collet,GHN,HNT}.
 
 In Section 2 we discuss the set-up, describe the martingale approach we use, prove some general results on variance and  discuss the SP property and Assumption C. Section  3 gives our results  under the assumption of quasi-compactness in H\"older norms and also some applications. In Section 4
 we give our results when we have quasi-compactness of the transfer operator  in the bounded variation norm, and we give applications to piecewise expanding maps in higher dimensions.  The last section is a concluding discussion, while the Appendices
 describe the Gal-Koksma lemma we use and show that  Assumption C is satisfied for generic points in many of our applications.

\section{The  setup.}

We suppose that $(T,X,\mu)$ is an ergodic dynamical system. 
Let the transfer operator $P$ be defined by $\int \phi \psi\circ T d\mu=\int P \phi \psi  d\mu$ for all $\phi,~\psi\in L^2 (\mu)$ so that $P$ is  the adjoint of
the Koopman operator $U\phi:=\phi\circ T$ with respect to the invariant measure $\mu$. Suppose $\mathcal{B}_{\alpha}$ is a Banach space of functions
and  $\|\phi\|_{1}\le C\|\phi\|_{\alpha}$ where $\|.\|_{\alpha}$ is  the Banach space norm and $\|.\|_1$ is the $L^1$ norm with respect to $\mu$.
We assume  $P$ restricts to an operator $P:\mathcal{B}_{\alpha}\to \mathcal{B}_{\alpha}$ such that  
$\|P^n \phi \|_{\alpha} \le C_1\theta^n \|\phi\|_{\alpha}$ for all $\phi \in \mathcal{B}_{\alpha}$ such that $\int\phi~d\mu=0$.
This implies   exponential decay of correlations  of the form, that for some $0<\theta <1$, 
\[
|\int \phi \psi\circ T^n~d\mu -(\int\phi~d\mu)(\int\psi~d\mu)|\le C\theta^n  \|\phi\|_{\alpha}\|\psi\|_{1}
\]
for  all $\phi\in \mathcal{B}_{\alpha}$, $\psi\in \mathscr{L}^1{\mu}$.  In our applications we will  have the pairs 
$(BV(X),\mathscr{L}^1(\mu))$ or $(H_{\gamma} (X), \mathscr{L}^{1} (X))$ where  $BV(X)$ is  the 
space of function of bounded variation and $H_{\gamma} (X)$ is the space of H\"older functions of exponent $\gamma$.  For example if $T$ is  a smooth  uniformly expanding map of the unit interval  $X$ then  $\mathcal{B}_{\alpha}$
 could be taken  as the Banach space of functions of bounded variation $BV(X)$.  In this paper
we will consider Lipschitz rather than H\"older functions, as our results and proofs immediately generalize to the H\"older setting with the
obvious changes.

\begin{rmk}
The weaker assumption of exponential decay of correlations
\[
|\int \phi \psi\circ T^n~d\mu -(\int\phi~d\mu)(\int\psi~d\mu)|\le C\theta^n  \|\phi\|_{\alpha}\|\psi\|_{\infty}
\]
implies that $\|P^n \phi \|_1 \le C\theta^n  \|\phi\|_{\alpha}$ (by taking $\psi$ to be sign($P^n\phi)$) and hence $P$ contracts exponentially in the $\mathscr{L}^1$ norm.
This assumption is sufficient for all our results on variance in Section 2, with the exception of the proof of the boundedness of the terms $w_j$, given in Lemma~\ref{w_j:bound} which 
seems to require our stronger assumption that $\|P^n \phi \|_{\alpha} \le C_1\theta^n \|\phi\|_{\alpha}$. These estimates on the growth of   $w_j$ are used in the proof of Theorem~\ref{Lip}. If $\mathcal{B}_{\alpha}$ is the space of functions of  bounded variation then the  $w_j$ terms are easily seen to be uniformly bounded under the assumption
$\|P^n \phi \|_{BV} \le C\theta^n  \|\phi\|_{BV}$.
\end{rmk}

Let $p\in X$ and let $B_n(p)$ be a sequence of nested balls
about $p$ such that $\mu (B_n(p))=\frac{1}{n}$. 
 Let $1_{B_n(p)}$ 
be the characteristic function of $B_n (p)$. We will sometimes write $E[\phi ]$ or $\int \phi$  for the integral $\int \phi~d\mu$
when the context is understood.  
Our results generalize immediately to sequences of nested balls with bounds $\frac{C_1}{n^{\gamma_1}}\le \mu (B_n(p))\le \frac{C_2}{n^{\gamma_2}}$
for constants $C_1,~C_2>0$ and $0<\gamma_2\le \gamma_1  \le 1$ (only the norming constants change)
but for simplicity we discuss in detail only the case $\mu (B_n(p))=\frac{1}{n}$.

$1_{B_n(p)}$ may not lie in $\mathcal{B}_{\alpha}$ but we assume we may take an approximation to it, $\tilde{\phi}^{\alpha}_n$ such that:
\begin{itemize}
\item[(i)] $|1_{B_n(p)} - \tilde{\phi}^{\alpha}_n|_1 \le  \frac{1}{n^3}$ and;
\item[(ii)]   $ \|\tilde{\phi}^{\alpha}_n\|_{\alpha} \le C n^k $ where $C,~k$ are independent of $n$;
\item[(iii)] $\tilde{\phi}^{\alpha}_n\ge 0$,~$\tilde{\phi}^{\alpha}_n\ge\tilde{\phi}^{\alpha}_{n+1}$
\end{itemize}

\begin{rmk}\label{remark:delta}
If we are taking a H\"older approximation then condition (ii) is satisfied for  the balls $B_i=B(p,r_i)$ if there exists $\delta (p)>0$  and $C>0$ such that 
$ \mu\{ x: r< d(x,p) < r+\epsilon \} < C\epsilon^{\delta (p)}$.  This condition is satisfied if the invariant measure $\mu$ has a density $h$ with respect to 
Lebesgue measure $m$  such that  $h\in\mathscr{L}^{1+\eta} (m)$ for some $\eta>0$.

\end{rmk}

 We define   $\phi^{\alpha}_n = \tilde{\phi}^{\alpha}_n-\int \tilde{\phi}^{\alpha}_n$ so that  $\int \phi^{\alpha}_n  =0$.
For ease of notation we will  subsequently drop the superscript $\alpha$ on $\phi^{\alpha}_n$ and $\tilde{\phi}^{\alpha}_n$.


Define $\phi_0=0$ and for $n\ge 1$
\[
w_n =  P \phi_{n-1}+P^2 \phi_{n-2}+\ldots +P^n \phi_0=\sum_{j=1}^nP^j\phi_{n-j}
\]
so that $w_1=P\phi_0$, $w_2=P\phi_1 +P^2\phi_0$, $w_3=P\phi_2+P^2\phi_1+P^3\phi_0$ etc...
For $n\ge 1$ define
\[
\psi_n= \phi_n - w_{n+1} \circ T + w_n
\]
Recall our assumptions $ \|\tilde{\phi}_n\|_{\alpha} \le C n^k $(so $\| \phi_n\|_\alpha \le \tilde{C} n^k$) and $\|P^n \phi \|_{\alpha} \le C_1\theta^n \|\phi\|_{\alpha}$ for all $\phi \in \mathcal{B}_{\alpha}$ such that $\int\phi d\mu=0$.
Hence $\|w_n\|_{\alpha}\le C_2 \|\phi_n\|_{\alpha}$ ,
$\|w_n\circ T\|_{\alpha}\le C_3 \|\phi_n\|_{\alpha}$ (since $\|UP\phi\|_{\alpha}\le C \|\phi \|_{\alpha}$ for all $\phi \in \mathcal{B}_{\alpha}$)
 and  hence $\|\psi_n \|_{\alpha}\le C_4 \|\phi_n\|_{\alpha}$.
Using the fact that $P(w_{n+1}\circ T)=w_{n+1}P1=w_{n+1}$ one may show  that $P\psi_n=0$.

Since $UP(\cdot)=E[\cdot |T^{-1} \mathcal{B}]$, $P\psi_j=0$ implies that $E[\psi_j | T^{-1}\mathcal{B}]=0$ and in turn
$E[\psi_j\circ T^{j} | T^{-1-j}\mathcal{B}]=0$ (since $T$ preserves $\mu$). Furthermore $\psi_j\circ T^j$ is $T^{-j}\mathcal{B}$  measurable
for all $j\ge 0$.

Following the approach of Gordin we will express $\sum_{j=1}^n \phi_j \circ T^j$ as the sum of a 
(non-stationary) martingale difference array and a controllable error term and then use the 
following Theorem 3.2 from Hall and Heyde~\cite{Hall_Heyde}:

\begin{thm}[Theorem 3.2~\cite{Hall_Heyde}]\label{Hall_Heyde}
Let $\{S_{n,i},\mathcal{F}_{n,i},1\le i\le k_n,n\ge1\}$ be a zero-mean square-integrable martingale difference array  with differences $X_{n,i}$ and let $\eta^2$ be
an almost sure finite random variable.  Suppose that:

(a) $\max_i |X_{n,i} | \rightarrow 0$ in probability;

(b) $\sum_{i} X_{n,i}^2 \rightarrow \eta^2$ in probability;

(c) $E(\max_i X_{n,i}^2)$ is bounded in $n$;

(d) the $\sigma$-fields are nested: $\mathcal{F}_{n,i}\subset \mathcal{F}_{n+1,i}$ for  $1\le i\le k_n$,~$n>1$.

Then $S_{n,k_n}\rightarrow Z$ (in distribution) where the random variable $Z$ has the characteristic function  $E(\exp (-\frac{1}{2} \eta^2  t^2))$. 
\end{thm}

As  is common  in the  application of martingale theory to non-invertible dynamical systems we will have to consider the natural extension
so that we have a martingale in backwards time. 
We  outline our scheme of proof.

Let $(\sigma,\Omega,m)$ be the natural extension of $(T,X,\mu)$. Each $\psi_j$ lifts to to  a function
$\psi_j^{*}$ on  $\Omega$ in a natural way, $\psi_j^{*}(\ldots \omega_{-2} \omega_{-1} .\omega_0 \omega_1 \ldots ):=\psi_j (\omega_0)$. To simplify
notation we write simply $\psi_j$ instead of $\psi_j^{*}$.

We define scaling constants 
by $a_n^2=E(\sum_{j=1}^{n}\phi_j \circ T^j )^2$. This sequence of constants play the role of non-stationary variance. 
Giving estimates on the growth and non-degeneracy of $a_n$ in this non-stationary setting is more difficult than in the usual stationary  case.

We define a triangular array $X_{n,i}=\frac{1}{a_n} \psi_{n-i}\circ \sigma^{-i}$,
$i=1,\ldots,n, n\in\mathbb{N}$, and put $S_{n,i}=\sum_{j=1}^{i}X_{n,j}$ for the partial sums 
(along rows). 
 Then $X_{n,i}$ is $\mathcal{F}_i :=\sigma^i \mathcal{B}_0$ measurable where
 $\mathcal{B}_0$ is the $\sigma$-algebra $\mathcal{B}$ lifted to $\Omega$. Note that in Theorem~\ref{Hall_Heyde} we take $\mathcal{F}_{n,i}:=\mathcal{F}_i$ for all $n$ and $k_n=n$.
The $\mathcal{F}_i$ form an increasing sequence of $\sigma$-algebras.
  We obtain $E[S_{n,i+1}|\mathcal{F}_i]=S_{n,i}+E[X_{n,i+1}|\mathcal{F}_i]$ where
   by stationarity $E[X_{n,i+1}|\mathcal{F}_i]=E[\psi_{n-i-1}|\sigma^{-1}\mathcal{B}_0]=0$. 
Hence $E[S_{n,i+1}|\mathcal{F}_i]=S_{n,i}$
and for every $n\in\mathbb{N}$  $X_{n,i}$ is a martingale difference array with respect to $\mathcal{F}_i$.

We will then verify conditions $(a)$, $(b)$, $(c)$ and $(d)$ of Theorem~\ref{Hall_Heyde}. 
The hard part  lies in establishing $(b)$. This is in contrast
with the stationary setting where condition $(b)$ is usually a straightforward consequence 
of the ergodic theorem. Condition $(b)$ is 
established in~\cite{NSV} by using ~\cite[Lemma 3.3.]{Sethuraman_Varadhan}, however
 in our setting the Lipschitz norms of the observations $\tilde{\phi_i}$
are unbounded and other techniques have to be used.  

Once we have established $(a)$, $(b)$, $(c)$ and $(d)$ it follows  that $\lim_{n\to \infty} \frac{1}{a_n} \sum_{j=0}^{n-1} \psi_j \circ T^j \to N(0,1)$ in distribution.
In the final step we show that $\frac{1}{a_n} \sum_{j=1}^n [w_j \circ T^j -w_j \circ T^{j+1}] \to 0$ in 
$\mathscr{L}^1$ which  implies that 
$\lim_{n\to \infty} \frac{1}{a_n} \sum_{j=0}^{n-1} \phi_j \circ T^j \to N(0,1)$ in distribution.

\subsection{Some lemmas on variance}

In this section we establish some preliminary results on the growth of the variance $E[(\sum_{j=1}^n \phi_j )^2]$ that will be useful in determining the 
scaling constants $a_n$.

For further reference let us notice that $\|P^n \phi \|_{\alpha} \le C_1\theta^n \|\phi\|_{\alpha}$ and $\|\phi\|_1\le C_1\|\phi\|_{\alpha}$ and that there exists a constant $a$ such that 
\begin{equation}\label{tail.estimate}
\|\sum_{j>a\log i} P^j \phi_i \|_1 \le \frac{1}{i^3}.
\end{equation}

\begin{lemma}\label{var}
\[
\limsup_{n\to\infty} \frac{1}{\log n} E(\sum_{i=1}^n \phi_i\circ T^i )^2 \ge 1
\]

\end{lemma}

\noindent {\bf Proof:}
By exponential decay of correlations and~(\ref{tail.estimate}) we get for the long term interactions:
\[
\sum_{j>a\log i +i}\left| \int \phi_i \circ T^i \phi_j \circ T^j \right| \le \frac{c_1}{i^2},
\] 
where we used exponential decay and our bound $\|\phi_j\|_1\le C_1 \|\phi_j\|_{\alpha}\le Cj^k$, where $C$,$k$ are independent of $j$. This bound is 
from assumption~(ii).
 Recall $\phi_j= \tilde{\phi}_j-\int \tilde{\phi}_j$
and  
$\|\tilde\phi_j\|_1 \le \frac{C_3}j$
(for some $C_3$). 
Thus for the short term interactions we get
$$
\sum_{j=i+1}^{i+a\log i}  \int \phi_i\circ T^i  \phi_j\circ T^j 
=\sum_{j=i+1}^{i+a\log i}\int\tilde{\phi}_i\circ T^i  \tilde{\phi}_j\circ T^j  +O(\frac{a\log i}{i^2})
$$
whence
\[
\sum_{i=1}^n \sum_{j>i} E[\phi_i\circ T^i\phi_j\circ T^j] 
=O(1) +\sum_{i=1}^n \sum_{j=i+1}^{i+a\log i}  E[\tilde{\phi}_i\circ T^i \tilde{\phi}_j\circ T^j].
\]
Since
\[
E(\sum_{i=1}^n \phi_i\circ T^i )^2
=\sum_{i=1}^n  E (\phi_i^2) +2 \sum_{i=1}^n \sum_{j>i} E[\phi_i\circ T^i\phi_j\circ T^j]
\]
and $\sum_{i=1}^n \sum_{j=i+1}^{i+a\log i}  E[\tilde{\phi}_i\circ T^i \tilde{\phi}_j\circ T^j] \ge 0$ the lemma is proved.
\qed

\begin{lemma}\label{lemma:a}
\[
\sum_{i=1}^n \sum_{j=i+1}^n \int  \phi_i \circ T^i \phi_j \circ T^j = \sum_{i=1}^n\int (\phi_iw_i)\circ T^i
\]
\end{lemma}

\noindent {\bf Proof:}
Recalling that $\phi_0=0$ this follows by a direct calculation and rearrangement of terms
 as
\begin{eqnarray*}
\sum_{i=1}^{n-1} \sum_{j=i+1}^n\int \phi_i \circ T^i \phi_j \circ T^j
&=&\sum_{j=2}^{n} \sum_{i=1}^{j-1}\int \phi_i \circ T^i \phi_j \circ T^j\\
&=&\sum_{j=2}^{n} \sum_{i=1}^{j-1}\int P^{j-i}\phi_i  \phi_j \\
&=&\sum_{j=2}^{n} \int(\sum_{i=1}^{j-1} P^{j-i}\phi_i)  \phi_j \\
&=&\sum_{j=2}^{n} \int w_j  \phi_j.
\end{eqnarray*}
\qed

\noindent The following lemma is the main result of this subsection: 
\begin{lemma}\label{lemma:c}
\[
a_n=E(\sum_{i=1}^n \phi_i\circ T^i )^2=\sum_{i=1}^nE[\psi_i^2] - \int w_1^2 +\int w_{n+1}^2
\]
\end{lemma}

\noindent {\bf Proof:} Let us first observe that factoring out yields
\begin{eqnarray*}
\psi_j^2&=&\phi_j^2+2\phi_j (w_j-w_{j+1}\circ T) +(w_j -w_{j+1}\circ T)^2\\
&=& \phi_j^2+2\phi_j (w_j-w_{j+1} \circ T) +w_j^2+w_{j+1}^2\circ T-2w_j w_{j+1}\circ T
\end{eqnarray*}
which when integrated leads to
\begin{eqnarray*}
\int \psi_j^2&=&\int \phi_j^2 +2 \int \phi_j (w_j-w_{j+1} \circ T) + \int w_j^2 +\int w_{j+1}^2-2\int w_j w_{j+1}\circ T\\
&=&\int \phi_j^2 + 2\int \phi_j w_j -2\int P\phi_j w_{j+1} + \int w_j^2 +\int w_{j+1}^2-2\int Pw_j w_{j+1}\\
&=&\int \phi_j^2 + 2\int \phi_j w_j -2\int P\phi_j w_{j+1} + \int w_j^2 +\int w_{j+1}^2-2\int (w_{j+1}-P\phi_j)w_{j+1}\\
&=&\int \phi_j^2 + 2\int \phi_j w_j +\int w_j^2 -\int w_{j+1}^2.
\end{eqnarray*}
Since by Lemma~\ref{lemma:a}
$$
a_n
=\sum_{i=1}^nE(\phi_i^2)+2\sum_{i=1}^n \sum_{j=i+1}^n \int  \phi_i \circ T^i \phi_j \circ T^j 
=\sum_{i=1}^n\left(E(\phi_i^2)+ 2\int (\phi_iw_i)\circ T^i\right)
$$
the statement follows by substituting $\int\psi_j^2-\int w_j^2+\int w_{j+1}^2$ for the terms inside
the sum on the RHS and then telescoping out the expected values of $w_j^2$.
\qed

\subsection{Property (SP)}

Several authors~\cite{Kleinbock_Margulis,Chernov_Kleinbock} have used a property derived from the Gal-Koksma theorem (see Appendix) to prove the SBC property for 
sequences of balls. Later we will show that in certain settings the (SP) property also implies a CLT.

Suppose $B_i$ are balls and let $f_i=1_{B_i}\circ T^i $. If
  \[
  \sum_{i=m}^n \sum_{j=i+1}^n E(f_i f_j)-E(f_i)E(f_j)  \le
  C \sum_{i=m}^n E(f_i)\hspace{1cm} (SP)
  \]
  for arbitrary integers $n>m$ then the balls are said to have the (SP) property.

\subsection{Short returns and Assumption~(C)}

In this section we discuss a condition on short return times first considered, to our knowledge, by P.~Collet~\cite{Collet}. We have called
it Assumption (C).  This condition has been used to 
establish extreme value statistics~\cite{Collet,HNT,GHN} and dynamical Borel-Cantelli lemmas~\cite{GNO,HNPV}. 

Suppose   $p\in X$ and $B_i (p)$ is a nested sequence of balls centered
at a point $p$, with $\lim_i \mu (B_i (p))=0$.

{\bf Assumption~(C):} We say $(B_i (p)) $ satisfies assumption~(C) if there
exists $\eta(p) \in(0,1)$ and $\kappa(p) >1$ such that for all $i$
sufficiently large
\[
\mu (B_i (p) \cap T^{-r} B_i (p)
)\le \mu(B_i (p))^{1+\eta}
\]
for all $r=1,\ldots, \log^\kappa i$.

If $(B_i (p)) $ satisfies assumption~(C) then we can say more about the behavior of the constants $a_n$.

\begin{lemma}\label{lemma:d}
Under Assumption~(C) there exists a constant $C_1$   so that
\[
\int \left|\phi_j w_j\right|\le \frac{C_1\log j}{j^{1+\eta}}.
\]

\end{lemma}

\noindent {\bf Proof:} By the contraction property of the transfer operator one has by~(\ref{tail.estimate})
for a sufficiently large constant $a$ 
\[
\sum_{i<j-a\log j}\int \phi_j P^{j-i} \phi_i \le \frac{1}{j^2}.
\]
Let $\phi_j=\tilde{\phi}_j- \int \tilde\phi_j $ where
$\tilde{\phi}_j$ is the $\mathcal{B}_{\alpha}$ approximation to $1_{B_j (p)}$ and note that 
$\|\tilde\phi_j\|_1, \|\phi_j\|_1\le \frac{c_1}j$
(for some $c_1$). 
Hence we obtain in the $\mathscr{L}^1$-norm:
(as $\tilde\phi_j\ge0$)
\begin{eqnarray*}
\int |\phi_j w_j|&\le&\sum_{n=1}^{a\log j}\left( \int \tilde{\phi}_j P^n \tilde{\phi}_{j-n}
+\int\tilde\phi_{j-n}^2+\int\tilde\phi_j\int P^n\tilde\phi_{j-n} +\int \tilde\phi_j\int \tilde\phi_{j-n}\right)
+\mathcal{O}\left(\frac1{j^2}\right)\\
&=&\sum_{n=1}^{a\log j} \left(\int \tilde{\phi}_j P^n \tilde{\phi}_{j-n} +3\mu(\tilde\phi_{j-n})^2\right)
+\mathcal{O}\left(\frac1{j^2}\right)\\
&=&\sum_{n=1}^{a\log j} \int \tilde{\phi}_j P^n \tilde{\phi}_{j-n}
+\mathcal{O}\left(\frac{\log j}{j^2}\right),
\end{eqnarray*}
where we used that $\int \tilde\phi_j =\mathcal{O}(j^{-1})$.
Now by assumption~(C) we have 
$$
\int  \tilde{\phi}_j P^n \tilde{\phi}_{j-n}\le \int \tilde{\phi}_{j-n}\circ T^n \tilde{\phi}_{j-n}
\le\mu(B_{j-n}\cap T^{-n}B_{j-n})\le \frac{c_2}{(j-n)^{1+\eta}},
$$
for $n\le a\log j$, and thus
$$
\sum_{n=1}^{a\log j} \int \tilde{\phi}_j P^n \tilde{\phi}_{j-n}\le  \frac{c_3}{j^{1+\eta}}a\log j,
$$
proving the lemma.
\qed

\begin{lemma}\label{lemma:e}
If $(B_i (p))$ satisfies Assumption~(C) then
\[
E(\sum_{i=1}^n \phi_i\circ T^i )^2=\sum_{i=1}^nE[\phi_i^2] +O(1)=\log n +O(1).
\]
\end{lemma} 
\noindent {\bf Proof:} Rearranging the sums yields by Lemma~\ref{lemma:a}
\begin{eqnarray*}
E(\sum_{i=1}^n \phi_i\circ T^i )^2
&=&\sum_{i=1}^n E[\phi_i^2] +2\sum_{i=1}^{n-1} \sum_{j=i+1}^n\int \phi_i \circ T^i \phi_j \circ T^j\\
&=& \sum_{i=1}^n E[\phi_i^2]  +2\sum_{j=2}^{n} \int w_j  \phi_j
\end{eqnarray*}
and hence the result follows by Lemma~\ref{lemma:d} as $\eta>1$ 
\qed

\subsection{Bounds on $w_j$}

We now assume that $\|\phi\|_{\infty} \le C  \|\phi\|_{\alpha}$ which under our assumption on the transfer operator  implies
that for a mean-zero function $\phi\in \mathcal{B}$, $\|P^n\phi\|_{\infty}\le C\theta^n \|\phi\|_{\alpha}$ for some $C$, $0<\theta<1$ independently of $\phi$. For example if $\|.\|_{\alpha}$ were the Banach space of H\"older functions of exponent $\gamma$ on the unit interval then 
$\|\phi\|_{\infty} \le C  \|\phi\|_{\alpha}$. In the BV or quasi-H\"older norm indicator functions are bounded, and the proof that $w_j$ is uniformly bounded 
is straightforward in this case. 

 \begin{lemma}\label{w_j:bound}
 Assume   $\|P^n\phi\|_{\infty}\le C\theta^n \|\phi\|_{\alpha}$ then there exist  a constant $C_2$
 such that  $  \|w_j\|_{\infty}<C_2$ for all $j$ where $E_n=\sum_{j=1}^n \mu (B_i)$.
 \end{lemma}
  
  \noindent {\bf Proof:}
For some   $a>0$ we can achieve 
$\sum_{j=a\log n}^n |P^j\phi_{n-j}|_\infty\le c_1\sum_{j=\lfloor a\log n\rfloor}^n\theta^j(n-j)^k
=\mathcal{O}(n^{-2})$
and in particular $|P^j\phi_{n-j}|_\infty=\mathcal{O}(n^{-2})$ for all $j\ge a\log n$ and all $n$.
As in the previous lemma let $\tilde\phi_j$ be the $\mathcal{B}_{\alpha}$ approximation for $1_{B_j}$
and $\phi_j=\tilde\phi_j-\mu(\tilde\phi_j)$.
 In view of the tail estimate it is only necessary to bound 
$\sum_{j=1}^{\lfloor a\log n\rfloor} P^j \phi_{n-j}$ independently of $n$.

 (i) Bound from below: Since $\phi_j \ge -\mu(\tilde\phi_j)=\mathcal{O}(j^{-1})$ one obtains
  $\sum_{j=1}^{\lfloor a\log n\rfloor} P^j \phi_{n-j}\ge\sum_{j=1}^{\lfloor a\log n\rfloor}  c_2(n-j)^{-\gamma_1} 
  \ge \frac{-c_3\log n}{n^{\gamma_1}}$ for some constants $c_2, c_3$ independent of $j$ and $n$. Hence $w_n\ge -c_4$ for some $c_4>0$ and all $n$
  (independent of $\gamma1, \gamma_2$).

(ii) Bound from above:  Since $1_{B_{j+1}}\le 1_{B_j}$ one has 
$\tilde\phi_{j+1}\le\tilde\phi_j$ and in particular $\mu(\tilde\phi_{j+1})\le\mu(\tilde\phi_{j})$.
Hence $\phi_{j+1}-\phi_{j}\le\mu(\tilde\phi_j)-\mu(\tilde\phi_{j+1})$ and (as $\phi_0=0$)
\begin{eqnarray*}
w_m-w_{m-1}&=&\sum_{j=1}^{m-1}P^j(\phi_{m-j}-\phi_{n-1-j})+P^m\phi_0 \\
&\le& \sum_{j=1}^{m-1}\left( \mu (\tilde{\phi}_{m-1-j})- \mu (\tilde{\phi}_{m-j})\right)\\
&\le&\sum_{j=1}^{\lfloor a\log m\rfloor} \left(\mu (\tilde{\phi}_{m-1-j})- \mu (\tilde{\phi}_{m-j})\right)+\mathcal{O}(m^{-2}).
\end{eqnarray*}
Consequently ($w_1=P\phi_0=0$)
\begin{eqnarray*}
w_n&=&\sum_{m=2}^n(w_m-w_{m-1})+w_1\\
&\le&\sum_{m=2}^n \left(\sum_{j=1}^{\lfloor a\log m\rfloor}\left(\mu (\tilde\phi_{m-1-j})
-\mu (\tilde\phi_{m-j})\right)
+\mathcal{O}(m^{-2})\right)\\
&=&\sum_{j=1}^{\lfloor a\log n\rfloor}\sum_{m=2\vee\lceil e^\frac{j}a\rceil}^n \left(\mu (\tilde\phi_{m-1-j})-\mu (\tilde\phi_{m-j})
+\mathcal{O}(m^{-2})\right)\\
&=&\sum_{j=1}^{\lfloor a\log n\rfloor}  \left(\mu (\tilde\phi_{2\vee\lceil e^\frac{j}a\rceil-j})-\mu (\tilde\phi_{n-j})
+\mathcal{O}((2\vee e^\frac{j}a)^{-1})\right)\\
&\le&C_3
\end{eqnarray*}
for a constant $C_3$ independent of $n$ because 
$$
\sum_{j=1}^{\lfloor a\log n\rfloor}\mu (\tilde\phi_{n-j})\le c_5\frac{a\log n}{n^{\gamma_2}}\to0
$$
as $n\to\infty$ and 
$$
\sum_{j=1}^{\lfloor a\log n\rfloor}  \mu (\tilde\phi_{2\vee\lceil e^\frac{j}a\rceil-j})
\le c_6\sum_{j=1}^{\lfloor a\log n\rfloor} (e^\frac{j}a)^{-\gamma_2}=\mathcal{O}(1)
$$
for  constants $c_5$, $c_6$ independent of $n$.
\qed

\section{Decay in Lipschitz versus $\mathscr{L}^{1}$}

We take  $\mathcal{B}_{\alpha}$ to be the space of Lipschitz functions, the arguments we give hold for H\"older norms with obvious modification.
We assume that the transfer operator $P$, when restricted to $Lip (X)$, contracts exponentially:
\begin{equation}\label{transfer.operator.lipschitz}
||P^n \phi||_{Lip} \le C \theta^n ||\phi||_{Lip}
\end{equation}
for all Lipschitz functions $\phi$   such that $\int\phi\, d\mu =0$ for some $\theta \in (0,1)$, where $\theta$ and $C$ independent of $\phi$.

This implies
\begin{equation}\label{correlation.lipschitz}
\left|\int \phi \psi\circ T^n d\mu  -E[\phi]E[\psi]\right|\le C \theta^n \|\phi\|_{Lip}~\|\psi\|_{L^{1}}
\end{equation}
for the same  $\theta \in (0,1)$ and $C$ independent of $\phi$,~$\psi$.

For a sequence of (nested) balls $B_i$ we 
put $E_n=\sum_{i=1}^n\mu(B_i)$ and $S_n=\sum_{i=1}^n1_{B_i}\circ T^i$ for the `hit counter'
for an orbit segment of length $n$. The sequence of balls $B_i$ satisfies the 
{\em strong Borel-Cantelli} (SBC) property if 
\begin{equation}\label{SBC}
\lim_{n\rightarrow\infty}\frac{S_n(x)}{E_n}=1
\end{equation}
for almost every $x\in X$.

\begin{thm}\label{Lip}
Assume  that the transfer operator, when restricted to $Lip (X)$, contracts exponentially
as in~(\ref{transfer.operator.lipschitz}) for some $\theta \in (0,1)$.

Suppose  $B_i (p) $ be nested balls  about a point $p$ with $\mu (B_i)=\frac{1}{i}$.
Let $a^2_n=E(\sum_{j=1}^n (1_{B_i} -\frac{1}{i}))^2$.\\
(I) If  the nested sequence of balls $(B_i(p))$ satisfies Assumption~(C) 
and the SBC property~(\ref{SBC}) then 
\[
a_n^2=\log n +\mathcal{O}(1)
\]
and
\[
\frac{1}{\sqrt{\log n}} \sum_{j=1}^n \left(1_{B_i} -\frac{1}{i}\right)\to N(0,1)
\]
in distribution.

(II) If   $(B_i(p))$ has the SP property then 
\[
\frac{1}{\sqrt{a_n}} \sum_{j=1}^n (1_{B_i}\circ T^i  -\frac{1}{i})\to N(0,1).
\]

\end{thm}
\vspace{2mm}

\noindent {\bf Proof:}
We will let $\phi_j=\tilde{\phi_j}- \int \tilde{\phi_j}$,where $\tilde{\phi_j}$ be a Lipschitz approximation to $1_{B_j}$, such that 
$$
\left\{\begin{array}{rcl}\|\tilde{\phi_j}-1_{B_j}\|_1&<&\frac{1}{j^2}\\
\|\tilde{\phi_j}\|_{Lip}&\le& Cj^k\\ \tilde{\phi}_j&\ge0\end{array}\right..
$$
We define $w_n=P\phi_{n-1}+P^2\phi_{n-2}+\ldots + P^n \phi_0$ and put
 $\psi_n= \phi_n -w_{n+1} \circ T + w_n$. Then
 $P\psi_n=P\phi_n-w_{n+1}+\sum_{j=2}^n P^j\phi_{n-j+1}=0$ which corresponds to 
 $\int \psi_n\chi\circ T\,d\mu=\int \chi P\psi_n\,d\mu=0$ for any integrable $\chi$.
 Note that  $||\phi_j||_{\infty} \le ||\phi_j||_{Lip}, ||\phi_j||_1 \le ||\phi_j||_{Lip}.$
 
 \begin{lemma}\label{lemma:w}
 There exist  constants $C_4,k,a$ so that\\
(I)  $\|w_n\|_{Lip}\le C_4 n^k$,\\
(II) $ \|w_n\|_{\infty} \le C_4 $,\\
(III) $ \|w_n\|_1\le C_4\frac{\log n}{n}$.
 \end{lemma}
 
 \noindent {\bf Proof of Lemma~\ref{lemma:w}.} (I) 
By the contraction of the transfer operator for Lipschitz continuous functions one obtains
$$
\|w_n\|_{Lip}\le\sum_{j=0}^\infty\|P^j\phi_n\|_{Lip}
\le \sum_{j=0}^\infty C_1\theta^j\|\phi_n\|_{Lip}\le c_1 n^k
$$

(II) Is a consequence of Lemma~\ref{w_j:bound}.

(III) For sufficiently large $a$  we get
\begin{eqnarray*}
||w_n||_1        &\le&  \sum_{j=1}^{a\log n} ||P^{j} \phi_{n-j}||_1+\sum_{j=a\log n+1}^{n} ||P^j \phi_{n-j}||_1\\
                      &\le& \sum_{j=1}^{a\log n} ||\phi_{n-j}||_1+\sum_{j=a\log n+1}^{\infty} ||P^j \phi_{n-j}||_{Lip}\\
                      &\le& \sum_{j=1}^{a\log n} ||\phi_{n-j}||_1+\sum_{j=a\log n+1}^{n} C_1 \theta^j ||\phi_{n-j}||_{Lip}\\
                      &\le& \frac{a\log n}{n-a\log n}+ c_4 \frac{\log^2 n}{n^2}\\
                      &\le&c_5\frac{\log n}n
\end{eqnarray*}
for some $c_4,c_5$ independent of $n$.

Now put $C_4=\max(c_1,c_5)$.
\qed

\vspace{3mm}


\noindent As before let $(\sigma,\Omega,m)$ be the natural extension of $(T,X,\mu)$ and
put  $a_n^2=E(\sum_{j=1}^{n}\phi_j \circ T^j )^2$ for the rescaling factors where the $\psi_j$ lift to 
$\Omega$ in a natural way.
By Assumption~(C), $a_n^2~\sim \log n$ by Lemma~\ref{lemma:e}.
Again we put $X_{n,i}=\frac{1}{a_n} \psi_{n-i}\circ \sigma^{-i}$,~$i=1,\ldots,n$  which are
  $\mathcal{F}_i =\sigma^i \mathcal{B}_0$  measurable where $\mathcal{B}_0$ is the $\sigma$-algebra 
  $\mathcal{B}$ lifted to $\Omega$. 
  The  $\mathcal{F}_i$ form an  increasing sequence of $\sigma$-algebras.
We put $S_{n,i}=\sum_{j=1}^{i}X_{n,j}$, $i=1,\dots,n$ ($k_n=n$), where the $X_{n,i}$
and obtain $E[S_{n,i+1}|\mathcal{F}_i]=S_{n,i}+E[X_{n,i+1}|\mathcal{F}_i]$ but by stationarity 
$E[X_{n,i+1}|\mathcal{F}_i]=E[\phi_{n-i-1}|\sigma^{-1}\mathcal{B}]=0$. 
Hence $E[S_{n,i+1}|\mathcal{F}_i]=S_{n,i}$
and   $X_{n,i}$ is a martingale difference array with respect to $\mathcal{F}_i$.

We now show condition (a) and (c) hold (clearly (d) holds). To see (a) and (c) calculate
$\int \psi_n^2 d\mu\le \|\psi_n\|_{\infty} \|\psi_n\|_{1}\le C\frac{\log^2n}{n}$.  Hence condition $(a)$ and $(c)$ hold.

We now prove (I) and  show that under Assumption~(C),  $\sum_{i=1}^nX_{n,i}^2\rightarrow 1$ in probability and  hence condition $(b)$ holds.

\begin{lemma}\label{peligrad}
$$
\frac1{\log n}\sum_{j=1}^n\psi_j^2\circ T^j\rightarrow 1
$$
in probability as $n\rightarrow\infty$.
\end{lemma}

\noindent {\bf Proof.}  We follow an argument given by Peligrad~\cite{Peligrad}.
As $\psi_j=\phi_j+w_j-w_{j+1}\circ T$ we obtain 
\begin{eqnarray*}
\psi_j^2 &=&\phi_j^2+2\phi_j w_j +w_j^2+w_{j+1}^2\circ T-2w_{j+1}\circ T(\phi_j +w_j)\\
&=& (\phi_j^2 +2\phi_j w_j +w_j^2+w_{j+1}^2\circ T-2w_{j+1}\circ T (\psi_j +w_{j+1}\circ T) \\
&=& \phi_j^2  +(w_j^2 -w_{j+1}^2 \circ T) -2\psi_j  w_{j+1} \circ T
+2\phi_j  w_j.
\end{eqnarray*}
We want to sum over $j=1,\dots, n$ and normalize by $\log n$ and wish to estimate
the error terms which are the last four terms on the RHS. The terms $w_j^2 -w_{j+1}^2 \circ T$ are bounded and telescope
so may be neglected. 


In order to estimate the third of the error terms, $\psi_j w_{j+1}\circ T$ we proceed
like Peligrad (page 9) using a truncation argument.
Let $w_j^{\epsilon}=w_j 1_{\{|w_j|\le \epsilon\sqrt{\log n}\}}$, where for simplicity of notation we have left out  the 
dependence on $n$. Then
$$
\int\left(\sum_{j=1}^n \psi_j \circ T^j w_{j+1}^{\epsilon}\circ T^{j+1}\right)^2
=\sum_{j=1}^n\int \left(\psi_j \circ T^j w_{j+1}^{\epsilon}\circ T^{j+1}\right)^2
\le \epsilon^2\sum_{j=1}^n\int \psi_j^2
$$
since the cross terms vanish (for $j>i$), as
\begin{eqnarray*}
 \int (\psi_jw_{j+1}^{\epsilon}\circ T)\circ T^{j} (\psi_iw_{i+1}^{\epsilon}\circ T)\circ T^{i} 
 &=& \int (\psi_jw_{j+1}^{\epsilon}\circ T)\circ T^{j-i} (\psi_iw_{i+1}^{\epsilon}\circ T)\\
 &=& \int (\psi_jw_j^{\epsilon}\circ T)\circ T^{j-i-1} P(\psi_iw_i^{\epsilon}\circ T)\\
&=& \int (\psi_jw_{j+1}^{\epsilon}\circ T)\circ T^{j-i-1} w_{i+1}^\epsilon P\psi_i=0
\end{eqnarray*}
as $P(\psi_i w_{i+1}^{\epsilon}\circ T)=w_{i+1}^{\epsilon}P\psi_i$.

For any $a>\epsilon$ we obtain using Tchebycheff's inequality (on the second term):
\begin{eqnarray*} 
P\left(\left|\frac1{\log n}\sum_{j=1}^n \psi_j\circ T^j w_{j+1} \circ T^{j+1}\right|>a\right)\hspace{-7cm}\\
&\le &P\left(\max_{1\le j\le n} \left|\frac1{\sqrt{\log n}}w_{j+1}\circ T^{j+1}\right|>\epsilon\right)+
P\left(\left|\frac1{\log n}\sum_{j=1}^n \psi_j \circ T^j w_{j+1}^{\epsilon}\circ T^{j+1}\right|>a\right)\\
&\le& P(\max_{1\le j\le n} |w_{j+1}\circ T^{j+1}|>\epsilon \sqrt{\log n})
+\frac{\epsilon^2}{a^2 \log n}\sum_{j=1}^n\int \psi_j^2\\
&=& P(\max_{1\le j\le n} |w_j\circ T^{j+1}|>\epsilon \log n)+c_1\frac{\epsilon^2}{a^2}.
\end{eqnarray*}
In the last line we used $\sum_{j=1}^n E [\psi_j^2] \sim \log n $ by  Lemma~\ref{lemma:c} and Lemma~\ref{lemma:e}.
By boundedness of the $w_j$ (Lemma~\ref{w_j:bound}) one gets
that $ P(\max_{1\le j\le n} |w_{j+1}\circ T^{j+1}|>\epsilon \sqrt{\log n} )\to 0$ 
for every $\epsilon>0$ as $n\rightarrow\infty$. Choosing $a=\epsilon^{\frac12}$ we 
conclude that 
$\frac1{a_n}\sum_{j=1}^n \psi_j\circ T^j w_{j+1} \circ T^{j+1}$ converges to zero in probability
as $n\rightarrow\infty$.

For the fourth error term $\frac1{{\log n}}2\sum_{j=1}^n(\phi_jw_j)\circ T^j$ we obtain by
 Lemma~\ref{lemma:e}:
$$
\left\|\sum_{j=1}^n(\phi_jw_j)\circ T^j\right\|_1
\le\sum_{j=1}^n\|\phi_jw_j\|_1
\le c_2\sum_{j=1}^\infty\frac{\log j}{j^{1+\eta}}<\infty
$$
uniformly in $n$. Thus $\frac2{{\log n}}\sum_{j=1}^n(\phi_jw_j)\circ T^j\rightarrow0$ in
probability.

Since the term $\frac1{{\log n} }\sum_{j=1}^n \phi_j^2 \circ T^j$ converges to $1$ almost surely by the SBC property
the proof is complete.
\qed

\vspace{2mm}

\noindent Lemma~\ref{peligrad} completes the proof of part~(I) of the theorem.
In order to show (II) we proceed as in the proof of~(I) except for  the verification of condition~(b).  
We will prove a SBC property for 
$\phi_j^2 +2w_j \phi_j$. Decomposing $\phi_j =\tilde{\phi}_j-\mu (\tilde{\phi}_j)$ and defining $\tilde{w}_j=P\tilde{\phi}_{j-1} +...+P^{[a\log j]}  \tilde{\phi}_{j-[a\log j]}$
we see that  $\| \phi_j^2-\tilde{\phi}_j^2\|_1\le \frac{C}{j^2}$ and $\|w_j\phi_j-\tilde{w}_j \tilde{\phi}_j\|_1\le \frac{C\log j}{j^2}$ so it suffices to consider
the sequence $\tilde{\phi}_j^2+2\tilde{w}_j \tilde{\phi}_j$. This is because $\sum_{j=1}^n E [\tilde{\phi}_j^2+2\tilde{w}_j \tilde{\phi}_j]=\sum_{j=1}^n E [\phi_j^2+2 w_j \phi_j]+O(1)$
and $\mu$ almost surely, $\sum_{j=1}^n (\tilde{\phi}_j^2\circ T^j +2(\tilde{w}_j \tilde{\phi}_j) \circ T^j )=\sum_{j=1}^n (\phi_j^2\circ T^j +2 (w_j \phi_j)\circ T^j)+O(1)$.

Note that both $\tilde{w}_j$ and $\tilde{\phi}_j$ are positive functions. 
Let $E_n:=\sum_{j=1}^n E [\tilde{\phi}_j^2+2\tilde{w}_j \tilde{\phi}_j]$.  We will use Proposition 8.1, a form of 
the Gal and Koksma theorem as stated by Sprindzuk (see Appendix) to show that $\frac{1}{E_n}\sum_{j=1}^n \tilde{\phi}_j^2\circ T^j +2(\tilde{w}_j \tilde{\phi}_j) \circ T^j \to 1$
almost surely. For this we want to use Proposition~\ref{prop:sprindzuk} with
 $f_j=\tilde\phi_j^2+2\tilde{w}_j\tilde\phi_j$, $g_j=\int f_j$ and $h_j$ to be determined below.
 We need to estimate the terms in $ \int\left(\sum_{i=m}^n\int \tilde\phi_j^2+2\tilde{w}_j\tilde\phi_j\right)^2$

 In order to verify the 
 condition of the proposition we look at the three individual sums as follows:\\
 (i) 
The fact that condition~(SP) holds for the functions $\tilde\phi_j$ implies
$$
\sum_{i=m}^n \sum_{j=i+1}^{n}\left| \int \tilde{\phi}_j \circ T^{j-i}  (\tilde{\phi}_i ) 
- E[\tilde{\phi}_j ]E[ \tilde{\phi}_i ]\right|
 \le C \sum_{i=m}^n E [\tilde{\phi}_j].
$$
Since $E(\tilde\phi_j^2)-E(\tilde\phi_j)=\mathcal{O}(j^{-k})$ we obtain
$$
\sum_{i=m}^n \sum_{j=i+1}^{n} \left| \int \tilde{\phi}_j^2 \circ T^{j-i}  \tilde{\phi}_i^2  
- E[\tilde{\phi}_j^2 ]E[ \tilde{\phi}_i^2]\right|
\le C \sum_{i=m}^n E [\tilde{\phi}_j] + \sum_{i=m}^n\mathcal{O}(i^{-k+1}).
$$
(ii) Lemma~\ref{lemma:w} ($|\tilde{w}_j|_\infty\le C_4\forall j$) now yields
\[
\sum_{i=m}^n \sum_{j=i+1}^{n}\left| \int (\tilde{\phi}_j  \tilde{w}_j)\circ T^{j-i}  (\tilde{\phi}_i  \tilde{w}_i)
- E[\tilde{\phi}_j \tilde{w}_j]E[ \tilde{\phi}_i \tilde{w}_i]\right|  
\le CC_4^2  \sum_{i=m}^n E [\tilde{\phi}_j] + \sum_{i=m}^n \mathcal{O}(i^{-k+1}).
\]
(iii) In the same way we obtain for the `mixed' terms
$$
\sum_{i=m}^n \sum_{j=i+1}^{n}\left| \int (\tilde{\phi}_j  \tilde{w}_j)\circ T^{j-i}  \tilde{\phi}_i^2
-   E[\tilde{\phi}_j\tilde{w}_j ]E[ \tilde{\phi}_i^2 ]\right|
   \le CC_4  \sum_{i=m}^n E [\tilde{\phi}_j] + \sum_{i=m}^n\left( \mathcal{O}(i^{-k+1})\right).
$$

Combining~(i), (ii) and~(iii) yields for all $m<n$ and some constant $c_1$:
$$
\int\left(\sum_{i=m}^n\int \tilde\phi_j^2+2\tilde{w}_j\tilde\phi_j\right)^2
\le c_1\sum_{i=m}^n\left(E(\tilde\phi_j)+\mathcal{O}(i^{-k+1})\right)
$$
which by Proposition~\ref{prop:sprindzuk} implies that 
$\frac{1}{E_n}\sum_{j=1}^n \left(\tilde{\phi}_j^2\circ T^j +2(\tilde{w}_j \tilde{\phi}_j) \circ T^j\right) \to 1$
almost surely, provided $k\ge2$.
\qed

\section{Applications to dynamical systems.}

Theorem~\ref{Lip} applies to a variety of dynamical systems including   Gibbs-Markov maps~\cite{Aaronson_Denker} 
and rational maps~\cite{Haydn}. For Gibbs-Markov maps
 it has been shown~\cite[Theorem 1]{GNO}
that nested sequences of balls $(B_i (p))$ satisfy both the Strong Borel Cantelli property and  assumption~C, so that (I) applies. 
 For rational maps~\cite[Theorem 10]{Haydn} shows that the transfer operator contracts exponentially in the $L^{\infty}$ norm hence if the (SP)
 property is also proved then (II) holds. More generally  (II) shows that 
proving the (SP) property for systems whose associated transfer operator has exponential decay suffices to prove the SBC property and the CLT for shrinking targets.

\section{Decay in $BV(X)$ versus $\mathscr{L}^1$}

It is known that summable decay of correlations in $BV(X)$ versus $\mathscr{L}^1$ implies the SP property by work of Kim ~\cite[Proof of Theorem 2.1]{Kim} (see also Gupta et al~\cite[Proposition 2.6]{GNO}).
Hence the statement in this setting is simpler.

Let the transfer operator $P$ be defined by 
$\int \phi \psi\circ T d\mu=\int P \phi \psi  \,d\mu$ for all $\phi,~\psi\in \mathscr{L}^2 (\mu)$,
that is  $P$ is  the adjoint of the Koopman operator $U\phi:=\phi\circ T$.

We assume  that the restriction of $P$ to the space $BV(X)$ is exponentially contracting, i.e.\
$P:BV(X)\to BV(X)$ satisfies
\begin{equation}\label{transfer.operator.BV}
\|P^n \phi \|_{BV} \le C\theta^n \|\phi\|_{BV}
\end{equation}
 for all $\phi \in BV(X)$ such that $\int\phi~d\mu=0$.

This implies  that $(T,X,\mu)$ has exponential decay of correlations  in $BV$ versus $L^1$, so  that for some $0<\theta <1$, 
\begin{equation}\label{correlation.BV}
\left|\int \phi \psi\circ T^n~d\mu -(\int\phi~d\mu)(\int\psi~d\mu)\right|\le C\theta^n  \|\phi\|_{BV}\|\psi\|_{1}
\end{equation}
for  all $\phi\in BV(X)$, $\psi\in \mathscr{L}^1 (\mu)$. In particular the measure $\mu$ is ergodic.

\begin{prop}\label{BV}
Assume the 
 transfer operator $P$ contracts exponentially as given by~(\ref{transfer.operator.BV})

Let  $B_i:=B(p,r_i)$ be nested balls of radius $r_i$ about a point $p$ such that $\mu (B_i)=\frac{1}{i}$,
and $a^2_n=E(\sum_{j=1}^n (1_{B_i} \circ T^i -\frac{1}{i}))^2$.
Then:\\
(I)   $\limsup_{n\to\infty} \frac{a_n}{\sqrt{\log n}}\ge 1$ and 
\[
\frac{1}{a_n} \sum_{j=1}^n (1_{B_i}\circ T^i  -\frac{1}{i})\to N(0,1).
\]\\
(II) If  the nested sequence of balls $(B_i(p))$ about $p$ satisfies Assumption~(C) then 
\[
a_n^2=E[
(\sum_{j=1}^n (1_{B_i}\circ T^i -\frac{1}{i}))^2]=\log n +\mathcal{O}(1)
\]
and
\[
\frac{1}{\sqrt{\log n}} \sum_{j=1}^n (1_{B_i}\circ T^i  -\frac{1}{i})\to N(0,1)
\]
in distribution.
\end{prop}

\noindent {\bf Proof:}
The proof is the same as for Theorem~\ref{Lip} with the simplification that the SP property holds automatically as we have  summable decay of correlations in $BV(X)$ versus $\mathscr{L}^1$ (see proof of~\cite[Theorem 2.1]{Kim}).
Furthermore  Lemma~\ref{var} shows that the variance is unbounded and Lemma~\ref{lemma:e} gives a precise rate of growth in the case
 that Assumption~(C) holds. 
\qed

\begin{rmk}
For one-dimensional maps of the interval, Proposition~\ref{BV} is  basically a consequence
of Conze and Raugi~\cite[Theorem 5.1]{Conze_Raugi}.  Follow the proof of~\cite[Theorem 5.1]{Conze_Raugi} taking $T_k=T$ for all $k$, $m$ to be the invariant measure $\mu$ and choosing
 $f_n=1_{B_n}(p)$. The rates of growth are given by  Lemma~\ref{var}  which
 shows that the variance is unbounded. Lemma~\ref{lemma:e} gives a precise rate of growth in the case that Assumption~(C) holds. In  Proposition~\ref{multidim} we extend
 these results to piecewise expanding maps in higher dimensions.

\end{rmk}

\section{Applications of Proposition~\ref{BV}.}

 Proposition~\ref{BV}  applies to certain classes of one-dimensional
maps such as  piecewise expanding  maps of the interval $T:X\to X$
with $\frac{1}{T'}$ of bounded variation and possessing an 
absolutely continuous invariant measure with density bounded away from zero (those maps satisfying the assumptions
of~\cite[Theorem 2.1]{Kim}, see also~\cite{GNO}). For these systems, Assumption~(C) has been 
shown to hold for nested balls about $\mu$ a.e.~$p\in X$~\cite{HNT,GHN}.
In the next subsection we generalize these results to piecewise expanding maps in higher dimensions.

\subsubsection{Piecewise expanding maps in higher dimensions}\label{subset:example-Saussol}

In this section we prove the Strong Borel Cantelli property and  the CLT for shrinking balls in a class
of expanding maps in higher dimensions. We also  show that assumption C holds for
$\mu$-a.e.\ point.

The Banach spaces will be given by  $\mathscr{L}^1$, defined with respect to the Lebesgue measure on $\mathbb{R}^n$, and a  quasi-H\"older space
with properties analogous to BV which we define below.  A key property of the quasi-H\"older space is that 
characteristic functions of balls have bounded norm (as in the BV norm) which turns out to be a very useful property. 

The maps are defined on compact sets $Z\in \R^N$. Denote by
$\dist(\cdot,\cdot)$  the usual metric in $\R^N$ and for $\eps>0$ let $B_{\eps}(x)=\{y\in\R^N: \dist(x,y)<\eps\}$
be the $\eps$-ball centred at $x$.  Let $B_\eps(A)=\{y\in \R^N: \dist(y,A)\leq \eps\}$ and
write $Z^\circ$ for the interior of $Z$ and $\overline{Z}$ its closure.

A map $T:Z\rightarrow Z$ is said to be a multidimensional piecewise expanding map, if there exists
 a  family of finitely many disjoint open sets $\{Z_i\}$ such that $\l(Z\setminus\bigcup_{i} Z_i)=0$ and there exist open sets $\widetilde{Z_i}\supset\overline{Z_i}$ and $C^{1+\alpha}$ maps $T_i: \widetilde{Z_i}\to\R^N$
 (for some  $0<\alpha\leq 1$) and some sufficiently small real number $\eps_1>0$  such that for all $i$,
\begin{itemize}
\item (H1) $T_i(\widetilde{Z_i})\supset B_{\eps_1}(T(Z_i))$ and $T_i|_{Z_i}=T |_{Z_i}$;
\item (H2) For $x,y\in T(Z_i)$ with $\dist(x,y)\leq\eps_1$,
$$|\det DT_i^{-1}(x)-\det DT_i^{-1}(y)|\leq c|\det DT_i^{-1}(x)|\dist(x,y)^\alpha;$$
\item (H3) There exists $s=s(f)<1$ such that $\forall x,y\in T(\widetilde{Z_i}) \textrm{ with } \dist(x,y)\leq\eps_1$, we have
$$\dist(T_i^{-1}x,T_i^{-1}y)\leq s\, \dist(x,y).$$
\item (H4) Let $G(\eps,\eps_1):=\sup_x G(x,\eps,\eps_1)$ where
\begin{equation}
G(x,\eps,\eps_1):=\sum_{i}\frac{\l(T_i^{-1}B_{\eps}(\partial T Z_i)\cup B_{(1-s)\eps_1}(x))}{\l(B_{(1-s)\eps_1}(x))}
\end{equation}
and assume that \begin{equation}
\label{sc}
\sup\limits_{\delta\leq\eps_1}\big(s^\alpha+2\sup\limits_{\eps\leq\delta}\frac{G(\eps)}{\eps^\alpha}\delta^\alpha\big)<1 \footnote{This condition could be greatly simplified as follows.  Suppose the boundaries of $Z_i$ are $C^1$ codimension one embedded compact submanifold, then define the quantity:
$$
\eta_0(T):=s^{\alpha}+\cfrac{4s}{1-s}Y(T)\cfrac{\gamma_{N-1}}{\gamma_{N}}
$$
where $$Y(T)=\sup_{x}\sum_i \# \left\{\text{smooth pieces
intersecting $\partial V_i$ containing $x$}\right\},$$ is the
maximal number of smooth components of the boundaries that can
meet in one point and $\gamma_N=\frac{\pi^{N/2}}{(N/2)!}$, the
$N$-volume of the $N$-dimensional unit ball of $\mathbb{R}^N$. We
require that $\eta_0(T)<1$, and this may replace the condition (\ref{sc}) above.}
\end{equation}
\end{itemize}

 We now introduce the Banach space of quasi-H\"older functions in which the spectrum of the  Perron-Frobenius operator $P$ is investigated.
Given a Borel set $\Gamma\subset Z$, we define the oscillation of $\varphi\in \mathscr{L}^1(\l)$ over $\Gamma$ as
$$\mathrm{osc}(\varphi,\Gamma):=\esssup\limits_{\Gamma}\varphi-\essinf\limits_{\Gamma}\varphi.$$

The function $x\mapsto \mathrm{osc}(\varphi,B_{\eps}(x))$ is measurable (see \cite[Proposition 3.1]{Ke})
For $0<\alpha\leq1$ and $\eps_0>0$, we define the  $\alpha$-seminorm of $\varphi$ as
$$
|\varphi|_{\alpha}=\displaystyle\sup_{0<\eps\leq\eps_0}\eps^{-\alpha}\int_{\R^\N}\mathrm{osc}(\varphi,B_{\eps}(x))\,\dif\l(x).
$$
Let us consider the space of functions with bounded $\alpha$-seminorm
$$
V_\alpha=\{\varphi\in \mathscr{L}^1(\l): |\varphi|_\alpha<\infty\},
$$
and endow $V_\alpha$ with the norm
$$
\|\cdot\|_\alpha=\|\cdot\|_1+|\cdot|_\alpha
$$
which makes it into a Banach space. We note that $V_\alpha$ is independent of the choice of $\eps_0$ and that $V_\alpha$ is continuously injected in $\mathscr{L}^{\infty}(\l)$.
According to \cite[Theorem 5.1]{SB}, there exists an absolutely continuous invariant probability measure
(a.c.i.p.) $\mu$, with density bounded above, and bounded below from zero,  which has exponential decay of correlations against $\mathscr{L}^1$ observables on
the finitely many mixing components of  $V_\alpha$: in view of the next Theorem \ref{piecewise} we will from now restrict ourselves to one of those components, by taking a mixing iterate of $T$.  More precisely, if the map $T$ is as defined above and if $\mu$ is the mixing a.c.i.p., then there exist constants $C<\infty$ and $\gamma<1$ such that
\begin{equation}\label{eq:doc-pen}
\Big|\int_{Z}\psi\circ T^n\, h\, \dif\mu -\int \psi \dif\mu \int h \dif \mu  \Big|\leq C\|\psi\|_{\mathscr{L}^1}\|h\|_\alpha\gamma^n
\end{equation}
for all $\psi\in\mathscr{L}^1$  and for all $h\in V_\alpha$. Moreover $\|P^n \phi \|_{\alpha} \le C \|\phi\|_{\alpha}$ for all
$\phi \in V_{\alpha}$ and thus equation~\ref{transfer.operator.BV} holds.

We now show that characteristic functions of balls are bounded in the $ \| \cdot \|_\alpha$ norm. 

\begin{lemma}\label{char_bound}
Let $B_i(p)$ be a nested sequence of balls about a point $p\in X$, then there exists a constant $C_3(\alpha)$ such that
\[
\| 1_{B_i} \|_\alpha \le C_3(\alpha)
\]
for all $i$.
\end{lemma}

\noindent{ \bf Proof:} 
Take any set $A$ with a rectifiable boundary. If $p$  is not in a  2$\epsilon$  neighborhood of the boundary of $A$, then the oscillation is zero, otherwise it is 1. Therefore we have
$\int osc(1_A, B_{\epsilon} (p) ) \,\dif\l(p)\le c_1 \epsilon$. Then we must divide by $\epsilon^{\alpha}$.
As $\alpha \le 1$ we have the ratio bounded by $c_1 * (\epsilon_0)^{1-\alpha}$.  
\qed

\vspace{2mm}

\noindent The boundedness of the characteristic functions in the $\|\cdot\|_\alpha$-norm allows
us to proceed as in Proposition~\ref{BV} (see also~\cite{Conze_Raugi}) and to 
obtain the following result.

\begin{prop}\label{multidim}
Assume a piecewise expanding map $T$ on a compact set $Z\subset \mathbb{R}^n$
satisfies conditions (H1)--(H4) and is mixing with respect to its absolutely continuous 
invariant measure $\mu$. 
Let  $B_i:=B(p,r_i)$ be nested balls of radius $r_i$ about a point $p$ such that $\mu (B_i)=\frac{1}{i}$.
Then the variance $a_n^2:=E[(\sum_{j=1}^n (1_{B_i}\circ T^i -\frac{1}{i}))^2]$
satisfies 
\[
\frac{a_n}{\sqrt{\log n}}\ge 1
\]
and
\[
\frac{1}{\sqrt{a_n}} \sum_{j=1}^n (1_{B_i}\circ T^i  -\frac{1}{i})\to N(0,1)
\]
in distribution.
\end{prop}

\noindent {\bf Proof:}
The SBC property (I) is immediate from the decay of correlations, Equation~\ref{eq:doc-pen} and the bound
$\| 1_{B_i} \|_\alpha \le C_3(\alpha)$ by the proof of Proposition~\ref{char_bound}.
The growth estimate follows from Lemma 2.3.
\qed

\vspace{2mm}

\noindent We now make an additional assumption. 
Suppose that we have $M$ domains of local injectivity for the map $T$; if we take the join
${\mathcal Z}^j:=\bigvee_{i=0}^{j-1}T^{-i}{\mathcal Z}$, where  ${\mathcal Z}$ denotes the
partition, mod-$0$, into the closed sets $\overline{Z_i}, i=1,\cdots,M$, then on each element
$Z^{(j)}_l, l=1,\cdots, |{\mathcal Z}^j|$, each of which is the closure of its interior, the map
$T^j$ is injective and of class $C^{1+\alpha}$ on an open neighborhood of  $Z^{(j)}_l$:
we call $\tilde{Z}^{(j)}_l$ such an  extension.  In order to prove condition (C) we require
a further assumption which is also called the {\em finite range structure}. We assume:
\begin{itemize}
\item (H5) Let
${\mathcal U}^{(j)}:=\{f^jZ_l^{(j)}, \forall l=1,\cdots, |{\mathcal Z}^j|\}$, and put ${\mathcal U}=\cup_{j=1}^{\infty}{\mathcal U}^{(j)}$. Then ${\mathcal U}$ consists of only finitely many subsets of $Z$ with positive Lebesgue measure, hence $U_m=\inf_{U\in\mathcal{U}}m(U)$ is bounded below.
\end{itemize}

\begin{lemma} Under the assumptions (H1)--(H5)  Assumption~(C) is satisfied.
\end{lemma}

\noindent {\bf Proof.} Denote
$$
E_k(\eps):=\{x;\dist(f^kx,x)\le \eps \}.
$$
By Lemma~\ref{lemma.condition.c}  (see Appendix) it is enough to prove that there exists $C>0$, $\delta>0$ such that for all $k$ and $\eps$,
$$
\mu(E_k(\eps)))<C\eps^{\tau}.
$$
We now fix $j$ and consider the cylinder, say,  $Z^{(j)}_l$. Let us suppose that $\{z_k\}_{k\ge1}$ is a sequence of points in $Z^{(j)}_l$ converging to $x\in Z^{(j)}_l$, namely $\text{dist}(z_k, x)\rightarrow 0$ when $k\rightarrow \infty$, and that $\text{dist}(T^j(z_k),x)\rightarrow 0$ for $k\rightarrow \infty$. With abuse of definition we  say that such a point $x$ is {\em fixed}.  If there are  points in the sequence $\{z_k\}_{k\ge1}$ which are on the boundary of $Z^{(j)}_l$, we think of $T^j$ as its $C^{1+\alpha}$ extension on
$\tilde{Z}^{(j)}_l$.  We want to show that in $\tilde{Z}^{(j)}_l$  there is only one fixed point $x$. By contradiction,
 suppose $y$ is another fixed point and $\{w_k\}_k$ a sequence converging to $y$ and whose $T^j$ images converge to $y$ as well. Suppose that $\tilde{Z}^{(j)}_l$ is a convex set in such a way the segment $[x,y]$ is contained in $\tilde{Z}^{(j)}_l$ \footnote{If not we could join $x$ and $y$ with a chain of segments contained each in $\tilde{Z^{(j)}_l}$: the argument will work again since the sum of the lengths of those segments is larger than the distance between $x$ and $y$ and this is what we need in bounding from below.}. We now fix $\eta$ small enough and take $k$ big enough and such that $\text{dist}(x,z_k), \ \text{dist}(x,T^j(z_k)), \  \text{dist}(y,w_k), \  \text{dist}(y,T^j(w_k))$,  are all smaller than $\eta$. We also put $D_{m,j}:=\inf\{||DT^j(x)||\}>1$, where the $\inf$ is taken over the points $x$ where the derivative is defined. The norm is the operator norm, which is strictly larger than $1$ since the map is uniformly expanding. Then we have
$$
\text{dist}(x,y)\ge \text{dist}(T^j(z_k),T^j(w_k))-\text{dist}(x,T^j(z_k))-\text{dist}(y,T^j(w_k))
$$
and by applying Taylor's formula
$$
\text{dist}(x,y)\ge D_{m,j}\text{dist}(z_k,w_k)-2\eta\ge D_{m,j}[\text{dist}(x,y)-2\eta]-2\eta
$$
which gives a contradiction, since $D_{m,j}>1$, by sending $\eta$ to $0$. Hence $x$ is the only fixed point.

Let us now take a measurable set $V\subset \tilde{Z}^{(j)}_l$ containing the fixed point  $x\in \tilde{Z}^{(j)}_l$. We require that the diameter of the image $T^j(V)$ be at most $\eps$; such an image will therefore be contained in the ball of center $T^j(x)$ and of radius $\eps$. The Lebesgue measure of this ball will be equal to $\gamma_N \eps^N$, where the factor $\gamma_N$ was defined in the preceding footnote. Then we have
$$
\l(B_\eps(x))=\gamma_N \eps^N\ge \l(T^j(V))\ge |\det(DT^j(\kappa))|\l(V)
$$
for a suitable point $\kappa\in\tilde{Z}_l^j$,
where in the last inequality we used a local change of variable and the continuity of $DT^j$, finally $\kappa$ is a point in $\tilde{Z}^{(j)}_l$. By distortion, we could replace this point by another one, say $\iota$ such that $\l(T^j({Z^{(j)}_l})=|\det(DT^j(\iota))|\l({Z^{(j)}_l})$. We therefore get (with the constant $B$ from~(BD))
$$
\l(V)\le \frac{\gamma_N \eps^N B}{|\det(DT^j(\iota))|}\le \frac{\gamma_N \eps^N B  \ \l({Z^{(j)}_l})}{U_m}
$$
Since the density of the absolutely continuous invariant measure $\mu$ is bounded from above (remember it is in $\mathscr{L}^{\infty}(\l)$), by, say, $h_M$, and since each $Z^{(j)}_l$ will contribute with at most one fixed point, by taking the sum over the $l$ we finally get
$$
\mu\{x; \text{dist}(T^jx,x)\}\le \frac{\gamma_N h_M\ \eps^N B  \ }{U_m}.
$$
and this bound is independent of $j$.
\qed

As a consequence of Lemma 2.8 we have,
\begin{thm}\label{piecewise}
Assume a piecewise expanding map $T$ on a compact set $Z\subset \mathbb{R}^n$
satisfies conditions (H1)--(H5) and is mixing with respect to its absolutely continuous 
invariant measure $\mu$. 
For $\mu$ a.e. $p$ if  $B_i (p) $ are nested balls about  $p$ such that $\mu (B_i)=\frac{1}{i}$.
Then
\[
a_n^2=E[
(\sum_{j=1}^n (1_{B_i}\circ T^i -\frac{1}{i}))^2]=\log n +\mathcal{O}(1)
\]
and
\[
\frac{1}{\sqrt{\log n}} \sum_{j=1}^n (1_{B_i}\circ T^i  -\frac{1}{i})\to N(0,1)
\]
in distribution.
\end{thm}

\section{Discussion.}

There are several natural questions remaining unanswered. In particular can the CLT for shrinking targets be proved for Anosov systems or non-uniformly hyperbolic
diffeomorphisms? Chernov and Kleinbock have proved the SBC property for balls in Anosov systems~\cite{Chernov_Kleinbock} but  the SBC property is unknown for non-uniformly hyperbolic diffeomorphisms.
More generally can a limit theory be developed for the statistics of non-stationary stochastic processes arising as observations (which change in time)
on deterministic dynamical systems which may also may evolve in time, such as sequential dynamical systems?

\section{Appendices}

\subsection{Gal-Koksma Theorem.}
We  recall the following result of Gal and Kuksma as formulated by  W. Schmidt~\cite{W1,W2} and stated
by Sprindzuk~\cite{Sprindzuk}:

\begin{prop}\label{prop:sprindzuk}
  Let $(\Omega,\mathcal{B},\mu)$ be a probability space and let $f_k
  (\omega) $, $(k=1,2,\ldots )$ be a sequence of non-negative $\mu$
  measurable functions and $g_k$, $h_k$ be sequences of real numbers
  such that $0\le g_k \le h_k \le 1$, $(k=1,2, \ldots,)$.  Suppose
  there exists $C>0$ such that
  \begin{equation} \label{eq:sprindzuk}
    \tag{$*$} \int \left(\sum_{m<k\le n}( f_k (\omega) - g_k)
    \right)^2\,d\mu \le C \sum_{m<k \le n} h_k
  \end{equation}
  for arbitrary integers $m <n$. Then for any $\epsilon>0$
  \[
  \sum_{1\le k \le n} f_k (\omega) =\sum_{1\le k\le n} g_k   +
  O (\Theta^{1/2} (n) \log^{3/2+\epsilon} \Theta (n)
  )
  \]
  for $\mu$ a.e.\ $\omega \in \Omega$, where $\Theta (n)=\sum_{1\le k
    \le n} h_k$.
\end{prop}

\subsection{Assumption~(C) for expanding systems}

In this appendix we show that if we  define 
\[
E_k (\epsilon):= \{x: d(T^k x,x)\le \epsilon\}
\]
and  if the invariant measure has a density bounded above with respect to Lebesgue then assumption C is valid.

\begin{lemma}\label{lemma.condition.c}
Suppose  $\mu$ has a density $\rho$  with respect to Lebesgue measure which satisfies $0<C_1<\rho <C_2$and  there exists $C>0$, $\delta>0$  such that  for all $k$, $\epsilon$,
\[
\mu (E_{k} (\epsilon))< C \epsilon^{\delta}
\] 
Then for $\mu$ a.e. $p\in X$ there
exists $\eta(p) \in(0,2)$ and $\kappa(p) >1$ such that for all $i$
sufficiently large
\[
\mu (B_i (p) \cap T^{-r} B_i (p)
)\le \mu(B_i (p))^{1+\eta}
\]
for all $r=1,\ldots, \log^\kappa i$.

\end{lemma}

\noindent {\bf Proof.} Let  $C_2=\max \rho (x)$, $C_1=\min\rho(x)$   where 
$\rho(x)=\frac{d\mu}{dm} (x)$ is the density of $\mu$  with respect to Lebesgue measure $m$.

Let $\sigma \ge 1$ and  $\gamma>\sigma$.
We choose $\epsilon_k$ so that  for all $x$  a ball of radius $\epsilon_k$ about $x$, denoted $B(x,\epsilon_k)$, satisfies
 $C_1/k^{\sigma} \le \mu (B(x,\epsilon_k) )\le C_2/k^{\sigma}$.

 Let $A_k:= \{ x:  d(T^j x,x)\le \epsilon_k \mbox{ for some } 1\le j\le \log (k)^{\rho}\}$.
 Evidently $A_k\subset \bigcup_{j=1}^{\log^\rho k}E_j$.
  By the estimate on $E_k (\epsilon)$ for all
 large $k$, $\mu (A_k) \le C \epsilon_k^{\tau}$  where $\tau<\delta$.
 Let
\[
F_k:=\{ x: \mu (B(x,\epsilon_k) \cap A_k) \ge  1/k^{\gamma}\}
\]
and define the Hardy-Littlewood maximal function  $M_k$ for $\phi(x)= 1_{A_k} (x)\rho(x)$ by 
\[
 M_k(x):=\sup_{a>0}\frac{1}{m(B_a(x))}\int_{B_a(x)} 1_{A_k}(y)\rho(y)\,dm(y).
\]
If $x\in F_k$ then $M_{k}> C_1 k^{\sigma-\gamma}$.

A theorem of Hardy and Littlewood (\cite{Folland} Theorem~3.17) states  that
\[
m( |M_k|>C)\le c_3 \frac{\|1_{A_k} \rho \|_1}{C}
\]
for some constant $c_3$, where $\|\cdot\|_1$ is the $\mathscr{L}^1$ norm with respect to $m$. 

Hence 
\begin{eqnarray*}
m(F_k)&\le& m(M_k > C_1 k^{\sigma -\gamma}) \\
&\le& \mu(A_{k}) C_1  k^{\gamma-\sigma} \\
&\le& k^{\gamma-\sigma(1+\tilde{\tau})}
\end{eqnarray*}
where $0<\tilde{\tau} < \tau$. We need to alter $\tau$ to $\tilde{\tau}$ to take into account the fact that 
a ball of radius $\epsilon$ has measure  roughly $\epsilon^{D}$.

Choosing $\sigma < \gamma < \sigma (1+\tilde{\tau})$ and $\sigma >1$ the series $\sum_k m(F_k)$ converges. 

So for $m$ a.e. $x_0 $ there exists an $N(x_0)$ such that $x_0\not \in F_k$ for all $k>N(x_0)$. 
Since $m(B(x,\epsilon_k)-m(B(x,\epsilon_{k+1}))\le \frac{2C_2}{k^2}$ this implies that for $\mu$ a.e. 
$x\in X$ there exists $\eta>0$ , $\kappa>0$  such that for all sufficiently large $i$, if $B_i(x)$ is a sequence of nested balls about $x$, $\mu (B_i (p))\sim\frac{1}{i}$ then   
\[
\mu (B_i(x) \cap T^{-r} B_i(x) )\le \mu (B_i (x))^{1+\eta}
\]
for $1\le  r\le \log (i)^{\kappa}$. This is Assumption~(C).
\qed

\section*{Acknowledgements}
MN and LZ were  partially
    supported by the National Science Foundation under Grant Number DMS-1101315. Part of this work was done while MN was visiting the {\em Centre de Physique Th\'eorique} in Marseille with a French CNRS support. SV was supported by the CNRS-PEPS {\em Mathematical Methods of Climate Theory}, by the ANR-Project {\em Perturbations} and by the {\sc PICS} ( Projet International de Coop\'eration Scientifique), {\em Propri\'et\'es statistiques des syst\`emes dynamiques det\'erministes et al\'eatoires}, with the University of Houston,  n. PICS05968. Part of this work was done while he was visiting the {\em Centro de Modelamiento Matem\'{a}tico, UMI2807}, in Santiago de Chile with a CNRS support (d\'el\'egation).


\begin{thebibliography}{00}

\bibitem{Aaronson_Denker} J. Aaronson, M. Denker {\em Local Limit Theorems for Gibbs-Markov Maps},{ Stoch. Dyn.} \textbf{1} no.2,  (2001),  193--237.

\bibitem{Boyarsky_Gora} A.~Boyarsky and P.~G\"ora. {\em Laws of Chaos: invariant measures and dynamical
systems in one dimension}. Birkhauser, 1997.


\bibitem{Chazottes_Collet} J.-R.~Chazottes and P.~Collet. {\em Poisson
  approximation for the number of visits to balls in non-uniformly
  hyperbolic dynamical systems}, {Preprint}.

\bibitem{Chernov_Kleinbock} N.~Chernov and D.~Kleinbock. {\em
  Dynamical Borel--Cantelli lemmas for Gibbs measures}, {Israel
  J. Math}. \textbf{122} (2001), 1--27.
  
  \bibitem{Collet} P. Collet. {\em Statistics of closest return for some
  non-uniformly hyperbolic systems}, {Erg. Th. \& Dyn. Syst.},
  \textbf{21} (2001), 401--420.
  
  \bibitem{Conze_Raugi}  Conze, J-P; Raugi, A. {\em Limit theorems for sequential expanding dynamical systems on [0,1]}
  Ergodic theory and related fields, \textbf{89-121}, Contemp. Math., 430, Amer. Math. Soc., Providence, RI, 2007. 

  
\bibitem{Dolgopyat} D.~Dolgopyat. {\em Limit theorems for partially
  hyperbolic systems}, Trans. AMS \textbf{356} (2004) 1637--1689.
  
\bibitem{Durret} R.~Durrett. {\em Probability: Theory and Examples},
  Second Edition, Duxbury Press, 2004.

\bibitem{Fayad} B.~Fayad, {\em Two remarks on the shrinking target
  property}, arXiv:math/0501205v1.
  
  \bibitem{Folland} G.~B.~Folland: {\it Real Analysis}: 2nd edition, Wiley 1999.
  
   \bibitem{Galatolo} S.~Galatolo, J.~Rousseau and B.~Saussol, {\em Skew products, quantitative recurrence,
  shrinking targets and decay of correlations},  Preprint.

\bibitem{Gordin69}
M.~I. Gordin. {The central limit theorem for stationary processes}.
  \emph{Soviet Math. Dokl.} \textbf{10} (1969) 1174--1176.

  
\bibitem{Gouezel} S.~Gou{\"e}zel, {\em A Borel--Cantelli lemma for
  intermittent interval maps}, {Nonlinearity}, \textbf{20} (2007),
  no.~6, 1491--1497.

\bibitem{GNO} C.~Gupta, M.~Nicol and W.~Ott, {\em A Borel--Cantelli
  lemma for non-uniformly expanding dynamical systems}, Nonlinearity
  23 (2010), no. 8, 1991--2008.
  
  \bibitem{GHN} C.~Gupta, M.~Holland and M.~Nicol, {\em Extreme value
  theory and return time statistics for dispersing billiard maps and
  flows, Lozi maps and Lorenz-like maps}, 
  {Erg. Th. Dyn. Syst.}, {\textbf 31-5}  (2011) 

  
  
  \bibitem{Hall_Heyde} P.~Hall and C.~C.~Heyde, {\em Martingale limit theory and its application}, {
  Probability and Mathematical Statistics}, Academic Press, 1980, New York.


\bibitem{Haydn} N.~Haydn, Haydn, {\em Convergence of the transfer operator for rational maps}, { Erg. Th. Dyn. Sys} \textbf{19} no.3,  (1999),  657Ð669.

\bibitem{HLV} N.~Haydn, Y.~Lacroix and S.~Vaienti, {\em Hitting and
  return time statistics in ergodic dynamical systems}, {
  Ann. Probab.}, \textbf{33}, (2005), 2043--2050.
  
  
\bibitem{HNPV} N.~Haydn, M.~Nicol, T.~Persson and S.~Vaienti. {\em A note on Borel--Cantelli lemmas for non-uniformly hyperbolic dynamical systems}, 
{Erg. Th. Dyn. Syst.}, to appear

\bibitem{Hirata} M.~Hirata, {\em Poisson Limit Law for Axiom-A
  diffeomorphisms}, {Erg. Thy. Dyn. Sys., 13 (1993)}, 533--556.

\bibitem{HNT} M. P. Holland, and M. Nicol and A. T\"or\"ok, {\em
  Extreme value distributions for non-uniformly hyperbolic dynamical
  systems}, {Trans. AMS}. {\textbf 364} (2012), no. 2, 661�688.

\bibitem{Katok_Strelcyn} A. Katok and J.-M. Strelcyn, {\em Invariant manifolds,
  entropy and billiards; smooth maps with singularities}, Springer
  Lecture Notes in Math.  {\bf 1222}, (1986).
  
  \bibitem{Kessebohmer} J.~Jaerisch, M.~Kesseb\"{o}hmer and Bernd O.~Stratmann. {\em A Fr\'ecehet law and an Erd\"{o}s-Philipp law for maximal cuspidal windings.}
  Preprint.
  
  \bibitem{Kleinbock_Margulis} D. Kleinbock and G. Margulis. {\em Logarithm laws for flows on homogeneous spaces}, { Inv. Math.} {\textbf 138}, (1999),
  451-494.

\bibitem{Kim} D.~Kim, {\em The dynamical Borel--Cantelli lemma for
  interval maps}, {Discrete Contin. Dyn. Syst.} \textbf{17} (2007),
  no.~4, 891--900.
  
  
\bibitem{Liverani96}
C.~Liverani. {Central limit theorem for deterministic systems}.
  \emph{{International Conference on Dynamical Systems}} (F.~Ledrappier,
  J.~Lewowicz and S.~Newhouse, eds.), Pitman Research Notes in Math.
  \textbf{362}, Longman Group Ltd, Harlow, 1996, pp.~56--75.





\bibitem{LSV} C.~Liverani, B.~Saussol and S.~Vaienti, {\em A
  probabilistic approach to intermittency}, {Ergodic Theory
  Dynam. Systems} \textbf{19} (1999), no.~3, 671--685.


  
\bibitem{Maucourant} F.~Maucourant, {\em Dynamical Borel--Cantelli
  lemma for hyperbolic spaces}, {Israel J. Math.}, \textbf{152}
  (2006), 143--155.
  
  \bibitem{NSV} P.~N\'andori, D.~Sz\'asz and T.~Varj\'u. {\em A central limit theorem for time-dependent dynamical systems}.
  Journal of Statistical Physics, DOI 10.1007/s10955-012-0451-8.

\bibitem{Peligrad}   M.~Peligrad, {\em Central limit theorem for triangular arrays of non-homogeneous Markov chains}
To appear in Prob. Theory and Related Fields.


\bibitem{Phillipp} W.~Phillipp, {\em Some metrical theorems in number
  theory}, {Pacific J. Math.} \textbf{20} (1967) 109--127.
  
  \bibitem{Rudin} W.~Rudin. {\em Real and Complex Analysis}, Third Edition, 1987, McGraw Hill.

\bibitem{Sprindzuk} Vladimir G. Sprindzuk, {\em Metric theory of
  Diophantine approximations}, V. H. Winston and Sons, Washington,
  D.C., 1979, Translated from the Russian and edited by Richard
  A. Silverman, With a foreword by Donald J. Newman, Scripta Series in
  Mathematics. MR MR548467 (80k:10048).
  
  \bibitem{W1} W.~Schmidt, {\em A metrical theory in diophantine approximation}, {Canad. J. Math}, \textbf{12},
  (1960), 619--631.
  
   \bibitem{W2} W.~Schmidt, {\em Metrical theorems  on fractional parts of sequences}, {Trans. Amer. Math. Soc.}, \textbf{110},
  (1964), 493--518.
  
  \bibitem{Sethuraman_Varadhan}  S.~Sethuraman and S.~R.~S Varadhan. {\em A martingale proof of Dobrushin's theorem
  for non-homogeneous Markov chains}, {Electronic journal of proabability}, \textbf{10}, (2005), 1221-1235.
 
 \bibitem{Viana}  M.~Viana. {Stochastic dynamics of deterministic systems}, Brazillian Math. Colloquium 1997, IMPA, Lecture Notes.
    
\end{thebibliography}
\end{document}